\documentclass[11pt]{amsart}
\usepackage{amscd, amsmath}
\usepackage{ amssymb}
\usepackage[margin=2.3cm]{geometry}

\newtheorem{theorem}{Theorem}[section]
\newtheorem{lemma}[theorem]{Lemma}
\theoremstyle{definition}

\newtheorem{proposition}[theorem]{Proposition}
\theoremstyle{remark}
\newtheorem{remark}[theorem]{Remark}
\numberwithin{equation}{subsection}
\theoremstyle{plain}

\newtheorem{corollary}{Corollary}
\newtheorem{question}{Question}

\def\F{\mathbb F}

\def\G{{\mathbb G}L}

\def\V{\mathbb V}
\def\W{\mathbb W}

\def\E{\mathbb E}
\def\F{\mathbb F}

\def\n{\mathfrak N}
\def\N{\mathcal N}

\def\V{\mathbb V}
\def\W{\mathbb W}

\copyrightinfo{2001}{enter name of copyright holder}
\newcommand{\secref}[1]{section~\ref{#1}}
\newcommand{\thmref}[1]{Theorem~\ref{#1}}
\newcommand{\lemref}[1]{Lemma~\ref{#1}}
\newcommand{\remref}[1]{Remark~\ref{#1}}
\newcommand{\propref}[1]{Proposition~\ref{#1}}
\newcommand{\corref}[1]{Corollary~\ref{#1}}
\newcommand{\eqnref}[1]{~{\textrm(\ref{#1})}}

\newcommand{\qref}[1]{Question~\ref{#1}}
\begin{document}
 \title[Existance of invariant non-degenerate bilinear form under a linear map]{On the existence of an invariant non-degenerate bilinear form under  a linear map}
\author[Krishnendu Gongopadhyay \and Ravi S. Kulkarni]{Krishnendu Gongopadhyay \and Ravi S. Kulkarni}
\address{Indian Institute of Science Education and Research (IISER) Mohali, Transit Campus: MGSIPAP Complex, Sector-26, Chandigarh 160019, India}
\email{krishnendu@iisermohali.ac.in, krishnendug@gmail.com}
\address{Department of Mathematics, Indian Institute of Technology Bombay, Powai, Mumbai 400076, India}
\email{punekulk@gmail.com}
\date{June 7, 2010}
\begin{abstract}
Let $\V$ be a vector space over a field $\F$. Assume that the characteristic of $\F$ is \emph{large}, i.e. $char(\F)>\dim \V$.  Let $T: \V \to \V$ be  an invertible linear map.  We answer the following question in this paper. 
\emph{When does $\V$ admit a $T$-invariant non-degenerate symmetric (resp. skew-symmetric) bilinear form?}  We also answer the infinitesimal version of this question. 

Following Feit-Zuckerman \cite{fz},  an element $g$ in a group $G$ is called \emph{real} if it is conjugate in $G$ to its own inverse. So it is important to characterize real elements in $\G(\V, \F)$. As a consequence of the answers to the above question, we offer a characterization of the real elements in $\G(V, \F)$.  

Suppose $\V$ is equipped with a non-degenerate symmetric (resp. skew-symmetric) 
bilinear form $B$.  Let $S$ be an element in the isometry group 
$I(\V, B)$. A non-degenerate $S$-invariant subspace $\W$ of $(\V, B)$ is called \emph{orthogonally indecomposable} with respect to $S$ if it is not an orthogonal sum of proper $S$-invariant subspaces. We classify the 
orthogonally indecomposable subspaces. This problem is nontrivial for the unipotent elements in $I(\V, B)$.  The \emph{level} of a unipotent  $T$  is the least integer $k$ such that $(T-I)^k=0$. We also classify the levels of  unipotents in  $I(\V, B)$. 
\end{abstract}

\maketitle

\section{Introduction}\label{intro}
Let $\F$ be a field, and let $\bar \F$ denote its algebraic closure.  
 Let $\V$ be a vector space of dimension $n$ over $\F$. The group of all invertible linear maps from $\V$ to $\V$ is denoted by $\G(\V, \F)$, or simply by $\G(\V)$ when there is no confusion about the field $\F$. In this paper we ask the following question.
\begin{question}\label{q1}
Given an invertible linear map $T: \V \to \V$, when does the vector space $\V$ admit a $T$-invariant non-degenerate symmetric (resp. skew-symmetric) bilinear form? 
\end{question}
Assuming that the characteristic of $\F$ is \emph{large}, i.e. $char(\F)>\dim \V$,  we have answered the question in this paper. 

Let $f(x)$ be a monic polynomial of degree $d$ over $\F$ such that $-1$, $0$, $1$ are not its roots. 
The \emph{dual} of $f(x)$ is defined to be the polynomial $f^{\ast}(x)=f(0)^{-1}x^d f(x^{-1})$. Note that if $f(x)=\Sigma_{i=1}^d a_i x^i$, then $f^{\ast}(x)=\frac{1}{a_0} \Sigma_{i=0}^d a_{d-i}x^i$. In other words, if $\alpha$ in $\bar \F$ is a root of $f(x)$ with multiplicity $k$, then $\alpha^{-1}$ is a root of $f^{\ast}(x)$ with the same multiplicity. The polynomial $f(x)$ is said to be \emph{self-dual} if $f(x)=f^{\ast}(x)$. Note that if $f$ is self-dual, then $a_0=1$. 

Let $T: \V \to \V$ be a linear transformation.  A $T$-invariant subspace is said to be \emph{indecomposable} with respect to $T$, or simply \emph{$T$-indecomposable} if it can not be expressed as a direct sum of two proper $T$-invariant subspaces. Clearly $\V$ can be written as a direct sum $\V=\Sigma_{i=1}^m \V_i$, where each $\V_i$ is $T$-indecomposable for $i=1,2,...,m$. In general, this decomposition is not canonical. But for each $i$, $(\V_i, T|_{\V_i})$ is ``dynamically equivalent'' to $(\F[x]/((p(x)^k), \mu_x)$, where $p(x)$ is an irreducible monic factor of the minimal polynomial of $T$, and $\mu_x$ is the operator 
$[u(x)] \mapsto [xu(x)]$.  Such $p(x)^k$ is an \emph{elementary divisor} of $T$. If $p(x)^k$ occurs $d$ times in the decomposition, we call $d$ the \emph{multiplicity} of the elementary divisor $p(x)^k$. By the theory of rational canonical form of linear maps, elementary divisors, counted with multiplicities, determine $(\V, T)$ upto ``dynamic equivalence'', cf. Kulkarni \cite{kulkarni} for the dynamical viewpoint. 

Let $\chi_T(x)$ denote the characteristic polynomial of an invertible linear map $T$. Let
 $$\chi_T(x)=(x-1)^e (x+1)^f \chi_{oT}(x),$$ where $e, f \geq 0$, and $\chi_{oT}(x)$ has no roots $1$, or $-1$. The polynomial $\chi_{oT}(x)$ is defined to be the \emph{reduced characteristic polynomial} of $T$. The vector space $\V$ has a $T$-invariant decomposition $\V=\V_1+ \V_{-1}+\V_o$, where for $\lambda=1, -1$, $\V_{\lambda}$ is the generalized eigenspace to $\lambda$, i.e.
$$\V_{\lambda}=\{v \in \V\;|\;(T-\lambda I)^n v=0\}, $$
and $T|_{\V_o}$ does not have any eigenvalue 1 or -1. Let $T_o$ denote the restriction of $T$ to $\V_o$. Clearly $T_o$ has the characteristic polynomial $\chi_{oT}(x)$.

\begin{theorem}\label{ansq1}
Let $\V$ be a vector space of dimension $n \geq 2$ over a field $\F$ of large characteristic. Let $T: \V \to \V$ be an invertible linear map.  Then $\V$ admits a $T$-invariant non-degenerate symmetric bilinear form if and only if the following  conditions hold.

(i) An elementary divisor of $T_o$ is either self-dual, or its dual is also an elementary divisor with the same multiplicity. 

(ii) Let $(x-1)^k$, resp. $(x+1)^k$, be an elementary divisor of $T$. Then either $k$ is odd or, if $k$ is even, then the multiplicity of the elementary divisor is an even number. So $2k \leq n$. $(\hbox{If $k$ is odd, then } k \leq n)$.
\end{theorem}

\begin{theorem}\label{ansq1'}
Let $\V$ be a vector space of dimension $2m \geq 2$ over a field $\F$ of large characteristic.  Let $T: \V \to \V$ be an invertible linear map.  Then $\V$ admits a $T$-invariant non-degenerate skew-symmetric bilinear form if and only if the following  conditions hold.

(i) An elementary divisor of $T_o$ is either self-dual, or its dual is also an elementary divisor with the same multiplicity. 
 
(ii) Let $(x-1)^k$, resp. $(x+1)^k$, be an elementary divisor of $T$. Then, either $k$ is even or,  if $k$ is odd, then the multiplicity of the elementary divisor is an even number. So $k \leq m$. $(\hbox{If $k$ is even, then } k \leq 2m)$.
\end{theorem}

When $T$  admits an invariant symmetric, resp. skew-symmetric bilinear form, then Part $(i)$ of the above theorems are implicit in the work on the conjugacy classes in orthogonal and symplectic groups, cf. Milnor \cite{milnor}, Springer-Steinberg \cite{ss}, Wall \cite{wall}. Part $(ii)$ of the theorems can be deduced from the work of Hesselink, cf. \cite[section-3]{hes}  but is not explicitly stated there.  The converse parts of the theorems are the really new contributions of our work. However, for completeness, we shall prove both parts in this paper, cf. \secref{pf}. As we shall see, the converse parts, i.e. the existence of invariant forms under conditions $(i)$ and $(ii)$,  require a subtle understanding of the arithmatic of field extensions and a detailed analysis of the unipotents.  

The analysis of the unipotents are mostly untouched in the work of  Milnor \cite{milnor}, Springer-Steinberg \cite{ss} or Wall \cite{wall}. Hesselink assumed the ground field to be ``quadratically closed" while dealing with the unipotents. While analyzing the unipotents, the connection with the Jacobson-Morozov lemma provides the required insight to  our work, cf. Gongopadhyay \cite{g}.  For the convenience of the reader, we provide here an {\it ab initio} complete treatment of this issue.   Notice that $sl(2, \F)$ admits a unique irreducible representation in each dimension, and further, this representation  admits a unique symmetric (resp. skew-symmetric) invariant form according as the dimension  is odd or even. In the terminology of physics, the so-called ``creation"- and ``annihilation"-operators are unipotent elements, exemplifying the terminal cases in the parts $(ii)$ in the above theorems.  On the other hand, an arbitrary unipotent element in the orthogonal or symplectic Lie algebra over an arbitrary field of large characheristic, by the Jacobson-Morozov lemma is contained in  some sub-algebra isomorphic to $sl(2, \F)$. This observation has motivated our precise formulation of \thmref{ansq1} and \thmref{ansq1'}.  Conversely, one could ask for an ``elementary" proof of the Jacobson-Morosov lemma based on \thmref{ansq1} and \thmref{ansq1'}. We hope to see such proof in the near future. 

\subsubsection*{Classification of real elements in the general linear group} Following Feit-Zuckerman \cite{fz}, an element $g$ in a group $G$ is said to be \emph{real},  if it is conjugate in $G$  to its own inverse.  If every element in the group $G$ is real, then $G$ is said to be \emph{real}. Reality properties of elements in linear algebraic groups are a topic of research interest due to their connection with the representation theory cf. Feit-Zuckerman  \cite{fz}, Moeglin et. al. \cite{mvw}, Singh-Thakur \cite{st1, st2}, Tiep-Zalesski \cite{tz}. It is an important problem to classify real elements in a group $G$. Wonenburger \cite{wonen} offered a characterization of real elements in $\G (n, \F)$ as a product of two involutions. As a corollary to \thmref{ansq1} and \thmref{ansq1'}, in the following we give a different criterion which classify real elements in $\G (n, \F)$. The corollary follows from the fact that if $T$ in $\G(n, \F)$ is real, then the characteristic polynomial $\chi_T(x)$ is self-dual. 

\begin{corollary}\label{cor1}
Let $\F$ be a field of large characteristic. Let $T$ be a real  element in $\G(n,\F)$.  Then $\V$ can be decomposed into a direct sum of $T$-invariant subspaces, $\V=\V_1 + \V_2$, such that $\V_1$ admits a $T$-invariant non-degenerate symmetric form and $\V_2$ admits a $T$-invariant non-degenerate skew-symmetric form. More precisely, 

(i) If $T$ has no eigenvalue $+1$ or $-1$, then $T$ preserves a non-degenerate symmetric, as well as a skew-symmetric bilinear form. 

(ii) Let $\chi_T(x)=(x-1)^n$, resp. $(x+1)^n$. If all the elementary divisors of $T$ are of even multiplicity, then it preserves a non-degenerate symmetric, as well as a skew-symmetric bilinear form. If all the elementary divisors are of odd multiplicity, then $T$ admits an invariant non-degenerate symmetric bilinear form and it can not admit any non-degenerate skew-symmetric form.

If some of the elementary divisors are of even multiplicity and some are of odd multiplicity, then there is a direct sum decomposition $\V=\V_1 + \V_2$ into $T$-invariant subspaces such that $\V_1$ admits a $T$-invariant non-degenerate symmetric form and $\V_2$ admits a $T$-invariant non-degenerate skew-symmetric form. 
\end{corollary}

We prove these theorems in \secref{pf}.  The existance of the invariant form under a unipotent map relies on a classification of the $T$-orthogonally indecomposable subspaces. 
\subsubsection*{Orthogonally indecomposable subspaces under unipotents}Let $\V$ be equipped with a non-degenerate symmetric or skew-symmetric bilinear form $B$.  The group of isometries of $(\V, B)$ is denoted by $I(\V, B)$. It is a linear algebraic group. When $B$ is symmetric, resp. skew-symmetric,  $(\V, B)$ is called a \emph{quadratic}, resp. \emph{symplectic} space. The group of isometries is denoted by $O(\V, B)$, resp. $Sp(\V, B)$. They are called the \emph{orthogonal} and the \emph{symplectic} groups respectively. A quadratic (resp. symplectic) space $(\V, b)$ is said to be a \emph{standard} quadratic (resp. symplectic) space if $\dim \V=2m$, and there exists subspaces $\W_1$ and $\W_2$ such that $\V=\W_1 + \W_2$, $\dim \W_1=\dim \W_2$ and  $b|_{\W_1}=0=b|_{\W_2}$. Here $\W_1$ and $\W_2$ are not unique, but $(\V, b)$ is unique upto isometry. The form $b$ is called a \emph{standard} symmetric, resp. skew-symmetric form. 

Let $\W$ be a non-degenerate $T$-invariant subspace of a quadratic (resp. symplectic) space. Then $\W$ is said to be  \emph{orthogonally indecomposable} with respect to an isometry $T$ if it is not an orthogonal sum of proper $T$-invariant subspaces. We classify the orthogonally indecomposable subspaces with respect to the unipotent isometries or their negatives.
\begin{theorem}\label{uniod}
Let $(\V, B)$ be a non-degenerate quadratic space over a field $\F$ of characteristic different from $2$. Let $T: \V \to \V$ be an isometry with minimal polynomial either $(x-1)^k$, or $(x+1)^k$. Let $\W$ be an orthogonally indecomposable subspace with respect to $T$. Then $\W$ has one of the following types.

(i) $\W$ is odd dimensional indecomposable. 

(ii) $\W$ is a standard space and each summand is indecomposable.
\end{theorem}
\begin{theorem}\label{uniod2}
 Let $(\V, B)$ be a non-degenerate symplectic space over a field $\F$ of characteristic different from $2$. Let $T: \V \to \V$ be an isometry with minimal polynomial either $(x-1)^k$, or $(x+1)^k$. Let $\W$ be an orthogonally indecomposable subspace with respect to $T$. Then $\W$ has one of the following types.

(i) $\W$ is even dimensional indecomposable. 

(ii) $\W$ is a standard space and each summand is indecomposable.
\end{theorem}
Note that in the case (i) of the above theorems, $\W$ is almost `standard' except for an one-dimensional summand! 
\begin{remark}
The tools we have used to derive \thmref{ansq1}---\thmref{uniod2} are simple and linear algebraic.  However,  in the large characteristics, \thmref{uniod} and  \thmref{uniod2}  can also be derived by a notable application of the Jacobson-Morozov lemma. For the statement of the Jascobson-Morozov lemma cf.  Bruhat \cite{bruhat}. In fact, our first proof of these theorems was based on the Jacobson-Morozov lemma. The condition \emph{large characteristic} on the base field implies the following. 

(a) For each self-dual elementary divisor $p(x)$ of $T_o$,  $p'(x)\neq 0$. Consequently $T_o$ has the Jordan-Wedderburn-Chevally decomposition, cf. Kulkarni \cite[Theorem-5.5]{kulkarni}. 

(b) When $T$ is unipotent,  the Jacobson-Morozov lemma holds for $T$, cf. Bruhat \cite{bruhat}. 

These are crucial ingredients in the derivation of the main theorems. Condition (a) is necessary for the proofs.  The Jacobson-Morozov lemma provides a sophisticated analysis for the unipotents, cf. Gongopadhyay \cite{g} when the field is algebraically closed.    However, the unipotents can also be analyzed without using the Jacobson-Morozov lemma. We present a simple linear algebraic approach to analyse them.  This approach is much simpler than the previous attempts, cf. for example, Hesselink \cite{hes}. Moreover,  it is valid over any field of characteristic different from $2$. 
\subsubsection*{Level of unipotents in orthogonal and symplectic groups}Recall that the \emph{level} of a unipotent $T$ in a linear algebraic group is the least integer $k$ for which $(T-I)^k=0$. The levels unipotents in a linear representation of a linear algebraic group $G$ is an important invariant of the representation.   
 
\begin{theorem}\label{level}
 Let $(\V, Q)$ be a non-degenerate quadratic space of dimension $\geq 3$  over a field $\F$ of  characteristic different from $2$. Let the maximal dimension of a subspace on which $Q=0$ is $l$. Let $k$ be the level of a unipotent isometry. Then $k$ will be one of the following.

either $(a)$ $k \leq l$,

or $(b)$ if $k>l$ and dimension of $\V$ is $2l$, then $k$ is odd and $k \leq 2l-1$.

or $(c)$ if $k>l$ and dimension of $\V$ is $\geq 2l+1$, 
 then $k$ is odd and $k \leq 2l+1$.
\end{theorem}
\begin{theorem}\label{level2}
 Let $(\V, Q)$ be a non-degenerate symplectic space  over a field $\F$ of characteristic different from $2$. Let the maximal dimension of a subspace on which $Q=0$ is $l$. Let $k$ be the level of a unipotent isometry. Then $k$ will be one of the following.

either $(a)$ $k \leq l$,

or $(b)$ if $k>l$ and dimension of $\V$ is $\geq 2l$, 
 then $k$ is even and $k \leq 2l$.
\end{theorem}
\end{remark}

\subsubsection*{Infinitesimal version of \qref{q1}}Further we ask the `infinitesimal' version of \qref{q1}. A  bilinear form $b$ is said to be \emph{infinitesimally invariant under a linear map} $S: \V \to \V$, or simply \emph{infinitesinally $S$-invariant} if for all $x$, $y$ in $\V$,
$$B(Sx, y)+B(x, Sy)=0. $$
The linear maps which preserve $B$ infinitesimally form an abelian group under addition, and this group is denoted by $\mathfrak I(\V, B)$. In fact, $\mathfrak I(\V, B)$ is the Lie algebra of the algebraic group $I(\V, B)$.

Let $f(x)$ be any monic polynomial of degree $d$ over $\F$ such that $f(0) \neq 0$. Let $f^-(x)=(-1)^df(-x)$. Then $f^-(x)$ is called the \emph{additive dual} polynomial to $f(x)$.  A monic polynomial $f(x)$ is called \emph{additively self-dual} if $f(x)= f(-x)$. Let $p(x)$ be an elementary divisor of $S$. If $p(x)$ is not a power of $x$, we call it a \emph{non-trivial elementary divisor} of $S$.  
\begin{theorem}\label{iinv}
 Let $\V$ be a vector space of dimension $n \geq 2$ over a  field $\F$ of large characteristic.  Let $S: \V \to \V$ be a linear map.   Then $\V$ admits an infinitesimally $S$-invariant  non-degenerate symmetric bilinear form if and only if the following  conditions hold.

(i) A non-trivial elementary divisor of $S$ is either additively self-dual, or its additive dual is also an elementary divisor with the same multiplicity.   

(ii) If $x^k$ is an elementary divisor of $S$ and $k$ is even, then the multiplicity of the elementary divisor is an even number. So $k \leq n$. $(\hbox{If $k$ is odd, then }k \leq n)$.
\end{theorem} 

\begin{theorem}\label{iinv1}
Let $\V$ be a vector space of dimension $2m$ over a  field $\F$ of  large characteristic. Let $S: \V \to \V$ be a linear map.  Then $\V$ admits an infinitesimally $S$-invariant non-degenerate  skew-symmetric bilinear form if and only if the following  conditions hold.

(i) A non-trivial elementary divisor of $S$ is either additively self-dual, or its additive dual is also an elementary divisor with the same multiplicity.    
 
(ii) If $x^k$ is an elementary divisor of $S$ and $k$ is odd, then the multiplicity of the  elementary divisor is an even number. So $k \leq m$. $(\hbox{If $k$ is even, then }k \leq 2m)$.
\end{theorem}
The proofs of the above theorems are analogous to the proof of \thmref{ansq1} and \thmref{ansq1'}.   
We omit the proofs. 

%%%%%%%%%%%%%%%%%%%
\section{Preliminaries} 
\subsection{The Standard Form}\label{standard}
Let $\W$ be a vector space over a field $\F$. Let $\W^{\ast}$ be the dual space to $\W$.  There is a canonical pairing $\beta: \W^{\ast} \times \W \to \F$ given by
$$\hbox{ for } w^{\ast}\in \W^{\ast},\;v \in \W,\;\;\beta(w^{\ast}, v)=w^{\ast}(v).$$
Moreover $\beta$ is non-degenerate, i.e. for each $w^{\ast}$ in $\W^{\ast}$, there is a $v$ in $\W$ such that $\beta(w^{\ast}, v) \neq 0$, and for each $v$ in $\W$, there exists  $w^{\ast}$ in $\W^{\ast}$ such that $\beta(w^{\ast}, v) \neq 0$. 

For $T$ in $\G(\W)$ and $v \in \W$, $w^{\ast}\in \W^{\ast}$, define  
 $(T \bullet w^{\ast})(v)=w^{\ast}(T^{-1}v)$.  
This defines an action of $\G(\W)$ on $\W^{\ast}$ from the left. 
 Further we have 
$$\beta(T \bullet w^{\ast},Tv)=(T \bullet w^{\ast})(Tv)=w^{\ast}(T^{-1}Tv)=w^{\ast}(v)=\beta(w^{\ast},v).$$
In this sense $T$ preserves the pairing $\beta$. 

Now consider the vector space 
$$\V=\W^{\ast}+ \W.$$ 
The pairing $\beta: \W^{\ast} \times \W \to \F$ can be extended canonically to a symmetric (resp. symplectic) form $b$ on $\V$ defined as follows.

(i) For $w\in \W$, $b(w, w)=0,$

(ii) For $w^{\ast} \in \W^{\ast}$, $b(w^{\ast}, w^{\ast})=0$,

(iii) For $w^{\ast} \in \W^{\ast}$ and $v \in \W$,          $b(w^{\ast},v)=w^{\ast}(v)=b(v,w^{\ast})$, (resp. $-b(v,w^{\ast})$ ).  
Since $\beta$ is non-degenerate, we see that $b$ is a non-degenerate symmetric (resp. skew-symmetric)  form. Moreover every invertible linear transformation $T: \W \to \W$ gives rise to an isometry as follows. 
 \begin{proposition}\label{spq2}
 There is a canonical embedding of $\G(\W)$ into $I(\V, b)$.
\end{proposition}
\begin{proof}
Let $T: \W \to \W$ be an invertible linear map. Define the linear map $h_T: \V \to \V$ as follows
$$h_T(u)=\left \{ \begin{array}{ll}
Tu & \hbox{ if } u \in \W\\
T  \bullet u & \hbox{if } u \in \W^{\ast}
\end{array} \right.$$
Now observe that for $v \in \W$, $w^{\ast} \in \W^{\ast}$,  
$$b(h_Tw^{\ast},h_Tv)=(h_Tw^{\ast})(h_Tv)=
(T \bullet w^{\ast})(Tv)=w^{\ast}(T^{-1}Tv)=w^{\ast}(v)=b(w^{\ast},v).$$
This shows that $h_T$ is an isometry. The correspondence $T \mapsto h_T$ gives the desired\\ embedding.
\end{proof}

Let $End(\W)$ denote the Lie algebra over $\F$ of all linear endomorphisms on $\W$. Then there is an action of $End(\W)$ on $\W^{\ast}$ as follows: for a linear map $S$, and for $w^{\ast}$ in $\W^{\ast}$, $v$ in $\W$, 
 $$ S \circ w^{\ast}(v)=w^{\ast}(-Sv).$$

Under this action $S$ infinitesimally preserves $\beta$ : 
\begin{eqnarray*}
\beta(S \circ w^{\ast}, v)+\beta( w^{\ast}, Sv)&=& S \circ w^{\ast}(v) + w^{\ast}(Sv)\\
&=& w^{\ast}(-Sv) + w^{\ast}(Sv)=0.
\end{eqnarray*}
Let $\V=\W^{\ast} + \W$. 
 The pairing $\beta$ can be extended to a symmetric, resp. skew-symmetric bilinear form $b$ on $\V$ by similar constructions as described above. Let $\mathfrak I(\V, b)$ denote the additive group of all linear maps on $\V$ which infinitesimally preserve $b$. It turns out that it is the Lie algebra of the algebraic group $I(\V, b)$.  
\begin{proposition}
There is a canoninal embedding of $End(\W)$ into $\mathfrak I(\V, b)$. 
\end{proposition}
\begin{proof}
Define the linear map $h_S: \V \to \V$ as follows
$$h_S(v)=\left \{ \begin{array}{ll}
S(v) & \hbox{ if } v \in \W,\\
S  \circ v & \hbox{if } v \in \W^{\ast}.
\end{array} \right.$$
Then $h_S$ preserves $b$ infinitesimally : for $w^{\ast}\in \W^{\ast}$, $v \in \W$, 
\begin{eqnarray*}
b(h_S w^{\ast}, v)+b(w^{\ast}, h_Sv) &=& b(S \circ w^{\ast}, v) + b(w^{\ast}, Sv)\\
&=& w^{\ast}(-Sv) + w^{\ast}(Sv)=0. 
\end{eqnarray*}
Then $S \mapsto h_S$ is the desired embedding. 
\end{proof}

\subsection{The norm and the trace of a field extension}
Let $\E$ be a finite extension of the field $\F$ of degree $[\E: \F]$. We denote the field extension by $\E/\F$.  For $\alpha$ in $\E$, the map $\hat \alpha:\E \to \E$ defined by $\hat \alpha ( e)=\alpha e$ is $\F$-linear. The \emph{trace of $\alpha$} from $\E$ to $\F$, denoted by $Tr_{\E/\F}(\alpha)$, is the trace of the $\F$-linear operator $\hat \alpha$.  The \emph{ norm of $\alpha$} from $\E$ to $\F$, denoted by $N_{\E/\F}(\alpha)$, is defined to be the determinant of $\hat \alpha$. The trace is an $\F$-linear map from $\E$ to $\E$, i.e. for all $e$, $e'$ in $\E$ and $a, b$ in $\F$, 
$$Tr_{\E/\F}(ae + be')=aTr_{\E/\F} (e)+b Tr_{\E/\F}( e').$$
The norm is a multiplicative map, i.e. for all $e$, $e'$ in $\E$, $N_{\E/\F}(ee')=N_{\E/\F}(e)N_{\E/\F}(e')$. Also for all $a$ in $\F$, $N_{\E/\F}(ae)=a^{[\E\;:\; \F]}N_{\E/\F}(e)$. The \emph{trace form} $\mathfrak t$ on $\E$ is defined by $\mathfrak t(e, e')=Tr_{\E/\F}(ee')$. The trace form is non-degenerate if and only if the extension $\E/\F$ is separable cf. Roman \cite[Theorem-8.2.2, p-204]{roman}.  

\section{Invariant form under a linear map}\label{lemmas}
\subsection{Correspondence between symmetric and skew-symmetric forms}\label{os}
Let $\V$ be a vector space over a field $\F$.  Suppose $T$ in  $\G(\V)$ is such that it has no eigenvalue $1$ or $-1$. Let $B$ be a $T$-invariant non-degenerate symmetric bilinear form on $\V$. Define a bilinear form $B^T$ on $\V$ as follows:
$$\hbox{For }u, v \hbox{ in } \V, \;\;B^T(u, v)=B((T-T^{-1})u, v).$$ 
Note that 
\begin{eqnarray*}
B^T(u, v)&=&B((T-T^{-1})u, v)\\
&=&B(Tu, v)-B(T^{-1}u, v)\\
&=&B(u, T^{-1}v)-B(u, Tv), \hbox{ since $T$ is an isometry}\\
&=&B(u, T^{-1}v-Tv)\\
&=&-B(u, (T-T^{-1})v)\\
&=&-B((T-T^{-1})v, u), \hbox{ since $B$ is symmetric}\\
&=&-B^T(v,u).
\end{eqnarray*}
Thus $B^T$ is a $T$-invariant  non-degenerate skew-symmetric form on $\V$.  Also it follows by the same construction that corresponding to each $T$-invariant skew-symmetric form, there is a canonical $T$-invariant symmetric form. We summarize this discussion in a proposition.
\begin{proposition}\label{qsf}
Let $T$ be an element in $\G(\V)$. If $T$ has no eigenvalue $ 1$ or $-1$, then there exists a $T$-invariant non-degenerate symmetric bilinear form on $\V$ if and only if there exists a $T$-invariant non-degenerate skew-symmetric form on $\V$. 
\end{proposition}
\subsection{Invariant form under a unipotent map}\label{cs}
\begin{lemma}\label{cu}
 Let $\V$ be an $n$-dimensional vector space of dimension $\geq 2$ over a field $\F$ of characteristic different from $2$. Let $T: \V \to \V$ be a unipotent linear map. Suppose  $\V$ is equipped with a $T$-invariant symmetric, resp. skew-symmetric bilinear form $B$.  Let $\V$ be indecomposable with respect to $T$. 

(i) Then $B$ is either non-degenerate or $B=0$. 

(ii) Suppose $B$ is non-degenerate. Then $n$ must be odd, resp. even. 
\end{lemma}
\begin{proof}
Let $T$ be an unipotent linear map. Suppose the minimal polynomial of $T$ is $m_T(x)=(x-1)^n$. Then without loss of generality we can assume that $T$ is of the form 
\begin{equation}\label{t}T=\begin{pmatrix} 1 & 0 & 0 & 0 & .... & 0\\ 1 & 1 & 0 & 0 & .... & 0 \\ 0 & 1 & 1 & 0 & ....& 0\\ &\ddots &  & \ddots    & \\ 0 & 0 & 0 & 0 & ....&0 \\ 0 & 0 & 0 & 0 & ....& 1 \end{pmatrix}.\end{equation}
Suppose $T$ preserves a non-degenerate bilinear form $B$. In matrix form, let $B=(a_{ij})$.
Since $T$ preserves $B$, hence $T^t B T =B$.  
This gives the following relations: For $1 \leq i \leq n-1$, 
\begin{equation}\label{1}
 a_{i+1,n}=0=a_{n, i+1},
\end{equation}
\begin{equation}\label{2'}
a_{i,j} + a_{i,j+1}+a_{i+1,j}+a_{i+1,j+1}=a_{i,j},
\end{equation}
\begin{equation}\label{2}
 i.e. \  a_{i,j+1}+a_{i+1,j}+a_{i+1,j+1}=0. \end{equation}
From the above two equations we have, for $0 \leq l \leq m-1$ and $l+1 \leq i \leq n-1$,
\begin{equation}\label{3}
a_{i,n-l}=0=a_{n-l,i}.
\end{equation}
This implies that $B$ is a triangular matrix of the form 
\begin{equation}\label{can}
B= \begin{pmatrix} a_{1,1} & a_{1,2} & a_{1,3} & a_{1,4} & .... & a_{1, n-2} & a_{1,n-1} & a_{1,n}\\
a_{2,1} & a_{2,2} & a_{2,3} & a_{2,4} & .... & a_{2,n-2} & a_{2,n-1} & 0\\
a_{3,1} & a_{3,2} & a_{3,3} & a_{3,4} & .... & a_{3,n-2} & 0 & 0\\
\vdots& \vdots&\vdots  &  \vdots& .... &\vdots  & \vdots & \vdots \\
a_{n-1,1} & a_{n-1,2} & 0 & 0 & ....& 0 & 0 & 0 \\
a_{n,1} & 0 & 0 & 0 & ....& 0& 0 & 0 \end{pmatrix},\end{equation} 
where, \begin{equation*} \label{e1} a_{i+1, j} + a_{i,  j+1} + a_{i+1, j+1}=0. \end{equation*} 
Using \eqnref{3} we have for $1 \leq l \leq n-1$, 
\begin{equation}\label{4} a_{l,n-l+1}=-a_{l+1,n-l}.\end{equation}
Using of the above equation with \eqnref{2} yields,  for $1 \leq l \leq n-1$,
\begin{equation}\label{6} a_{l,n-l+1}=a_{l,n-l}+a_{l+1,n-l-1} \end{equation}

\subsection{Proof of (i)} Suppose $B$ is degenerate. Hence the determinant of the matrix $B$ must be zero. Without loss of generality,  in the form \eqnref{can} of $B$ we assume for $1 \leq l \leq n$, $a_{l,n-l+1}=0$. Now, this implies from \eqnref{6},  for $1 \leq l \leq n-1$, $a_{l, n-l}=-a_{l+1, n-l-1}$. Thus for $1 \leq l \leq n-1$
\begin{equation}\label{7} a_{l,n-l}=(-1)^{n-2l} \ a_{n-l,l} \end{equation} 
Suppose $B$ is symmetric, resp. skew symmetric.  Then \eqnref{7} implies that $n$ must be even, resp. odd. This is a contradiction to (i) above. Hence we must have for $1 \leq l \leq n-1$, $a_{l, n-l}=0$.  Continuing the process, we have $a_{i,j}=0$ for all $(i,  j) \neq (1,1)$.  Choose a basis $\{e_1,....,e_n\}$ of $\V$ such that $T$ and $B$ has the above forms with respect to the basis. Thus $B(e_1, e_1)=a_{1,1}$, and $T(e_1)=e_1$. The complementary subspace of $\F e_1$ is the radical of $B$ and is $T$-invariant. This contradicts the assumption that $\V$ is $T$-indecomposable. Hence we must have $a_{1,1}=0$. Hence $B=0$. 
\subsection{Proof of (ii)} Suppose $B$ is symmetric, resp. skew-symmetric and non-degenerate. 
From \eqnref{4} we have
  $$a_{l,n-l+1}=(-1)^{n+1-2l} \ a_{n-l+1,l}.$$
Thus we must have $n$ is odd, resp. even. 

 This completes the proof. 
\end{proof} 
\begin{lemma}\label{ind}
Let $\V$ be an $n$-dimensional vector space of dimension $\geq 2$ over a field $\F$ of characteristic different from $2$. Let $T: \V \to \V$ be a unipotent linear map. Let $\V$ be $T$-indecomposable. Then there exists a $T$-invariant non-degenerate symmetric, resp. skew-symmetric bilinear form on $\V$ if and only if $n$ is odd, resp. even. 
\end{lemma}
\begin{proof}
Without loss of generality, as in the proof of the previous theorem, assume $T$ is of the form \eqnref{t}. Then any bilinear form $B$ of the form \eqnref{can} is preserved by $T$.  Consequently, as in the above lemma,  $B$ is symmetric, resp. skew-symmetric, if and only if $n$ is odd, resp. even. 
\end{proof} 
\begin{remark}\label{r1}
In \eqnref{t} if we replace $1$ by a $k \times k$ identity matrix $I$, then the same procedure produces a $T$-invariant non-degenerate symmetric, resp. skew-symmetric bilinear form according as $\frac{n}{k}$ is odd, resp. even where $n=\dim \V$. In this case, in \eqnref{can}, each $a_{i,j}$  is replaced by a $k \times k$ matrix $A_{i,j}$ and they satisfy \eqnref{2} i.e. 
$$A_{i+1, j} + A_{i,  j+1} + A_{i+1, j+1}=0.$$
\end{remark}

\subsection{The Induced Form}
\begin{lemma}\label{p-inv}
Let $\V$ be a vector space over a field $\F$ of large characteristic.  
Let \hbox{$T: \V \to \V$} be a linear map with characteristic polynomial $\chi_T(x)=p(x)^d$, where $p(x)$ is irreducible over $\F$ and is self-dual.  Let $\V$ be indecomposable with respect to $T$. Then $\dim \V$ is even, and there exists a  $T$-invariant 
non-degenerate symmetric, as well as skew-symmetric bilinear  form on $\V$.
\end{lemma}
\begin{proof} 
Since $\V$ is $T$-indecomposable, $(\V, T)$ is dynamically equivalent to the pair 
\hbox{($\F[x]/((p(x)^d)$, $\mu_x$)}, where $\mu_x$ is the operator 
$\mu_: [u(x)] \mapsto [xu(x)]$, cf. Kulkarni \cite{kulkarni}. Hence without loss of generality we assume $\V=\F[x]/((p(x)^d)$, $T=\mu_x$. 

Suppose that the degree of $p(x)$ is $2m$. Let $y=x + \frac{1}{x}$. Then $x^{-m}p(x)$ is a polynomial with indeterminate $y$ over $\F$. We denote this polynomial in $y$ by $q(y)$. Since $p(x)$ is irreducible, $q(y)$ is also irreducible. Since the characteristic of $\F$ is large, note that $p'(x) \neq 0$. 

Let $\E=\F[x]/(p(x))$, and $\E_1=\F[y]/(q(y))$. Clearly $\E_1$ may be taken as a subfield of $\E$. As a field extension $\E$ has degree $2$ over $\E_1$. Since $p'(x)\neq0$,  we see that each of the extensions $\E/\F$, $\E/\E_1$ and $\E_1/\F$,  is separable. Note that $N_{\E/\E_1}$ defines a $T$-invariant $\E_1$-valued quadratic form on $\E$. We denote the corresponding symmetric bilinear form  by $\N_{\E/\E_1}$. Define $B: \E \times \E \to \F$ by 
$$B(\alpha, \beta)=Tr_{\E_1/\F}( (\N_{\E/\E_1}(\alpha,1)(\N_{\E/\E_1}(\beta,1)).$$
 Then $B$ defines a non-degenerate $T$-invariant symmetric bilinear form on $\E$. Correspondingly, there exists a non-degenerate $T$-invariant skew-symmetric form $\sigma$ on $\E$. 
This proves the theorem for $d=1$, i.e. when $T$ is semisimple.

Assume $d$ is at least $2$.  It follows from Kulkarni \cite{kulkarni} that there exists a basis such that 
$$T=\begin{pmatrix} M & O &  &  & &\\
I & M & & &  &\\
O & I & M & & &\\
& & I & M & & &  \\
& & & \ddots& \ddots & & \\
& & & &I & M
    \end{pmatrix},$$
where $M$ denote a $2m \times 2m$ matrix with characteristic polynomial $p(x)$ and $I$ denote the $2m \times 2m$ identity matrix. All other entries in the above matrix are zeros. The matrix $M$ is unique up to conjugacy.  Also, the Jordan-Wedderburn-Chevalley decomposition exists. Hence $T=T_s T_u$,   where $T_s$ is semisimple, $T_u$ is unipotent,  $T_sT_u=T_uT_s$, and such a decomposition is unique. So, we can assume $T_u$ is of the form 

\begin{equation}T_u=\begin{pmatrix} I & O & O & O & .... & O\\ I & I & O & O & .... & O \\ O & I & I & O & ....& O\\ &\ddots &  & \ddots    & \\ O &O & O & O & ....&O \\ O & O & O & O & ....& I \end{pmatrix}.\end{equation}

where each $I$ is the $2m \times 2m$ identity matrix, and $n=2md$.  As mentioned in \remref{r1}, there exists a non-degenerate symmetric, resp. skew-symmetric bilinear form $B_u$ according to $d=\frac{\dim \V}{\deg p(x)}$ is odd, resp. even. Further $B_u$ is of the form
\begin{equation}\label{can2}
B_u= \begin{pmatrix} A_{1,1} & A_{1,2} & A_{1,3} & A_{1,4} & .... & A_{1, d-2} & A_{1,d-1} & A_{1,d}\\
A_{2,1} & A_{2,2} & A_{2,3} & A_{2,4} & .... & A_{2,d-2} & A_{2,d-1} & O\\
A_{3,1} & A_{3,2} & A_{3,3} & A_{3,4} & .... & A_{3,d-2} & O & O\\
\vdots& \vdots&\vdots  &  \vdots& .... &\vdots  & \vdots & \vdots \\
A_{d-1,1} & A_{d-1,2} & O& O & ....& O & O & O \\
A_{d,1} & O & O & O & ....& O& O & O \end{pmatrix},\end{equation} 
where each $A_{i,j}$ is a $2m \times 2m$ matrix and 
\begin{equation}\label{m1}  A_{i+1, j} + A_{i,  j+1} + A_{i+1, j+1}=0. \end{equation} 
As we have seen in the previous section, $B_u$ can be made symmetric, resp. skew-symmetric if and only if $d$ is odd, resp. even. 

In the Jordan decomposition, let 
\begin{equation}
T_s=\begin{pmatrix} M & O & ...& O & O \\ O & M & ...& O & O \\ & & \ddots & &  \\O & O & ...& M & O \\ O & O & ...& O & M\end{pmatrix},
\end{equation}
where $M$ is a semisimple matrix with characteristic polynomial $p(x)$.  By the semisimple case above, there exists a $M$-invariant non-degenerate symmetric, as well as a skew-symmetric bilinear form. Let $b$ be a non-degenerate $M$-invariant symmetric bilinear form. Then in matrix representation $b$ is a $2m \times 2m$ symmetric non-singular matrix. 

Case (i):  Suppose $d$ is odd. Let $d=2k+1$. We want to construct $B_u$ so that it is non-degenerate, symmetric and is invariant under $T_s$. 

Assume, for all $i, \ j$, $A_{i, j}=A_{j, i}$. Consider an element $\alpha_{k,l}$ in $A_{i,j}$. Then $\alpha_{k,l}$ is the $(2m(i-1)+k, 2m(j-1+l)$-th entry of $B_u$. Then $\alpha_{l,k}$ in $A_{j,i}$ is the $(2m(j-1)+1, 2m(i-1)+k)$-th entry of $B_u$. If $\alpha_{l,k}=\alpha_{k,l}$, i.e. if $A_{i,j}$ is chosen to be symmetric, then $B_u$ is symmetric. Thus, for all $i, \ j$, if we choose $A_{i,j}$ to be symmetric, then $B_u$ would be symmetric.  

For $1 \leq l \leq d-1$, choose $A_{l, d-l+1}=(-1)^{l-1} b$. Note that $A_{d-l+1,l}=(-1)^{d-l}b=(-1)^{2k-(l-1)}b=(-1)^{l-1}b=A_{l,d-l+1}$. Hence the choice satisfies the assumed symmetricity of $B_u$. 

Using the above choice, we have from \eqnref{m1} the following set of equations:
\begin{eqnarray}
& &  A_{1,d-1} + A_{2, d-2}=b\\
& &  A_{2,d-2} + A_{3,d-3}=-b \\
& &  \vdots  \\
& &  A_{l,d-l} + A_{l+1,d-l-1}=(-1)^{l+1}b\\
& &  \vdots  \\
& &  A_{k-1,k+2}+A_{k, k+1}=(-1)^k b\\
& & A_{k, k+1} + A_{k+1,k}=(-1)^{k+1} b. 
\end{eqnarray}
By the assumption that $B_u$ is symmetric, we rewrite the last equation as $2 A_{k, k+1}=(-1)^{k+1} b$, i.e. $A_{k, k+1}=(-1)^{k+1} \frac{b}{2}$. Now by back-substitution we have,  for $1  \leq m \leq k-1$, 
$A_{k-m, k+m+1}=(-1)^{k-m+1} \  \frac{2m+1}{2}b$.  Next we have the following set of equations:
\begin{eqnarray}
& & A_{1,d-2}+A_{2,d-3}=\frac{2k-3}{2}b\\
& & A_{2,d-3}+A_{3,d-4}=-\frac{2k-5}{2}b\\
& & \vdots \\
& & A_{k-l+1, k+l+1}+A_{k-l+1,k+l}=(-1)^{k-l+2} \ \frac{2l+1}{2} b \\
& & \vdots \\
& & A_{k+1,k-1}+A_{k+2,k}=-A_{k+2,k-1}=(-1)^k \frac{3}{2} b.
\end{eqnarray}
We already have the initial value $A_{k, k+2}=A_{k+2,k}=(-1)^{k-1}b$. From the last equation we have
$A_{k+1,k-1}=A_{k-1,k+1}=(-1)^k \frac{5}{2} b$. Now by back substitution as in the previous case, we get other solutions. Continuing the process, we have for all $i, j$, $(i, j) \neq (1,1)$, $A_{i, j}$ as a multiple of $b$ by a rational number over $\F$. Put $A_{1,1}=b$. Then $B_u$ is symmetric and non-degenerate. Now, each of the $2m \times 2m$ block of $B_u$ is invariant under $M$. Hence $T_s B_u T_s^t=B_u$. Thus $B_u$ is invariant under $T_s$, as well as it is invariant under $T_u$. Hence $B_u$ is a non-degenerate symmetric bilinear form invariant under $T$. 

Thus when $d$ is odd, $T$ has a invariant non-degenerate symmetric bilinear form. Consequently by \propref{qsf}, there exists a non-degenerate invariant skew-symmetric bilinear form. 

Case (ii):  Suppose $d$ is even, let $d=2k$. In this case $B_u$ can not be symmetric. We wish to construct a non-degenerate skew-symmetric bilinear form $B_u$ which is invariant under both $T_s$ and $T_u$. 

Assume for all $i$, $A_{i,i}=O$, and for all $i, \ j$, $A_{i, j}=-A_{j,i}$. 
As in the above case, we see that for all $i, \ j$, if $A_{i, j}$ is chosen to be symmetric, then $B_u$ would be skew-symmetric. Hence we choose for $1 \leq l \leq k$, $A_{l, n-l+1}=(-1)^{l-1}b$. Now using \eqnref{m1} and following similar procedure as above, we obtain a skew-symmetric non-degenerate bilinear form $B_u$ which is invariant under both $T_u$ and $T_s$, and hence it is also $T$-invariant. Consequently, by \propref{qsf}, there also exists a $T$-invariant non-degenerate symmetric bilinear form. 

This completes the proof. 
\end{proof}

The following proposition follows immediately from the above lemma. 
\begin{proposition}\label{p-inv1}
Let $\V$ be a vector space over a field $\F$ of large characteristic. Let $T: \V \to \V$ be a linear map with characteristic
polynomial $\chi_T(x)=\Pi_{i=1}^k p_i(x)^{d_i}$, where for each $i=1,2,...,k$, $p_i(x)$ is self-dual and is irreducible over $\F$. Then there exists a non-degenerate $T$-invariant symmetric, as well  skew-symmetric  bilinear form on $\V$.
\end{proposition}

\subsection{Proof of \thmref{uniod}}
For a unipotent isometry $T: \V \to \V$ with minimal polynomial $(x-1)^k$, we observe that $-T: \V \to \V$ is also an isometry with minimal polynomial $(x+1)^k$, and also the converse holds. Hence it is enough to prove the theorem for unipotents. 

Let $B$ be symmetric.  Let $T:\V \to \V$ be a unipotent isometry in $I(\V, B)$. Let $\W$ be a  $T$-indecomposable and orthogonally indecomposable subspace of $\V$.  Since $B|_{\W}$ is non-degenerate, we see that $\dim \W$ must be odd by \lemref{ind}. 

Suppose $\dim \W$ is even. Then $B|_{\W}=0$ by \lemref{cu}. Since $B$ is non-degenerate, there is a $T$-invariant $T$-indecomposable subspace $\W'$ such that $B|_{\W'}=0$, $\dim \W=\dim \W'$ and $\W + \W'$ is non-degenerate. Hence the multiplicity of the  elementary divisor $(x-1)^k$ must be even. 
This completes the proof of \thmref{uniod}.

\subsection{Proof of \thmref{uniod2}}
The proof is similar as above. 

\section{Proofs of \thmref{ansq1} and \thmref{ansq1'}}\label{pf}
\begin{lemma}\label{lem1}Suppose $T$ admits an invariant non-degenerate symmetric (resp. skew-symmetric) bilinear form $B$. Suppose $T$ has no eigenvalue $1$ or $-1$. Then the minimal polynomial $m_T(x)$ of $T$ is self-dual. 
\end{lemma}
\begin{proof}
 Note that for $u$, $v$ in $\V$ we have $B(Tu, v)=B(u, T^{-1}v)$. Using this identity it follows that 
for any $f(x)$ in $\F[x]$, $B(f(T)v, w)= B(v, f(T^{-1})w)$. Applying this to the minimal polynomial of $T$ we have for all $v$ in $\V$,  $m_T(T^{-1})v=0$. Thus if $\lambda$ in $\bar \F$ is a root of the minimal polynomial, then $\lambda^{-1}$ is also a root. Thus the minimal polynomial of $T$ is self-dual. 
\end{proof}

Let $p(x)$ be an irreducible factor of the minimal polynomial of $T_o$. \lemref{lem1} implies that $p(x)$ is either self-dual, or there is an irreducible factor $p^{\ast}(x)$ such that $p(x)$ is dual to $p^{\ast}(x)$. Suppose $p(x) \neq p^{\ast}(x)$. Let $\V_p=\hbox{ ker } p(T)^n$, $\V_{p^{\ast}}=\hbox{ ker } p^{\ast}(T)^n$, where $n=\dim \V$. Then $B|_{\V_p}=0=B|_{\V_{p^{\ast}}}$ and $B|_{\V_p+\V_{p^{\ast}}}$ is non-degenerate.
\subsection{Proofs of \thmref{ansq1} and \thmref{ansq1'}}
Suppose $T$ admits an invariant non-degenerate symmetric (resp. skew-symmetric) bilinear form $B$.  
  Let $\oplus$ denote the orthogonal sum and $+$ denote the usual sum of subspaces. It can be seen easily that there is a \emph{primary decomposition} of $\V$ (with respect to $T$) into $T$-invariant non-degenerate subspaces:
\begin{equation}\label{decom}
\V=\oplus_{i=1}^{k_1} \V_i \bigoplus \oplus_{j=1}^{k_2} \V_j
\end{equation}
where for $i=1,2,...,k_1$, $p_i(x)$ is self-dual, $\V_i=\V_{p_i}$, and $B|_{\V_i}$ is non-degenerate; for $j=1,2,...,k_2$, $\V_j=\V_{p_j} + \V_{p^{\ast}_j}$, $B|_{\V_{p_j}}=0=B|_{\V_{p_j^{\ast}}}$, here $p_j(x) \neq p_j^{\ast}(x)$.  

 So, without loss of generality assume $\V$ is of the form $\V_i$ or $\V_j$, i.e. $m_{T_o}(x)$ is of either of the form $p(x)^d$ or $q(x)^d q^{\ast}(x)^d$, where $p(x)$, $q(x)$, $q^{\ast}(x)$  are irreducible over $\F$, and $q^{\ast}(x)$ is dual to $q(x)$.

{\it Case (i).} Suppose $m_{T_o}(x)=p(x)^d$, where $p(x)$ is self-dual and irreducible over $\F$. By the theory of Jordan-canonical form, we have a direct sum decomposition: $\V=\Sigma_{i=1}^k \V_{d_i}$, where  $1 \leq d_1 <....< d_k=d$, and for each $i=1,...,k$, $p(x)^{d_i}$ is an elementary divisor, $\V_{d_i}$ is $T$-invariant and is free over the algebra $\F[x]/(p(x)^{d_i})$. We claim that each $\V_{d_i}$ is non-degenerate. For this it is sufficient to show that $\V_d$ is non-degenerate. The non-degeneracy of the other summands will follow by induction.  

If possible, suppose $\V_d$ is degenerate. Let $R(\V_d)$ be the  radical of $\V_d$, that is, 
$$R(\V_d)=\{v \in \V_d\;|\; B(v, \V_d)=0\}.$$
Let $v$ be a non-zero element in $R(\V_d)$. Since $R(\V_d)$ is $T$-invariant, let $p(T)v=0$. By the theory of elementary divisors it follows that there exist a $u$ in $\V_d$ such that $p(T)^{d-1}u=v$. Then for all $i<d$, and $w$ in $\V_{d_i}$, 
$$B(p(T)^{d-1}u, w)=B(u, p(T^{-1})^{d-1}w)=0.$$
Hence $v$ is orthogonal to $\V$, a contradiction to the non-degeneracy of $B$. Thus $\V_d$ must be non-degenerate. Thus for each $i$, $\V_{d_i}$ is non-degenerate. Note that the minimal polynomial of $T|_{\V_{d_i}}$ is $p(x)^{d_i}$. Hence $p(x)^{d_i}$ is self-dual. 

{\it Case (ii).}  $m_{T_o}(x)=p(x)^dp^{\ast}(x)^d$, where $p(x)$, $p^{\ast}(x)$ are irreducible over $\F$, $p(x) \neq p^{\ast}(x)$ and $\V=\V_p + \V_{p^{\ast}}$, $\dim \V_p=\dim \V_{p^{\ast}}$, 
$B|_{\V_p}=0=B|_{\V_{p^{\ast}}}$. By the theory of Jordan Canonical form, $\V$ has a direct sum decomposition:
$\V=\Sigma_{i=1}^k \V_{p^{d_i}}+ \Sigma_{j=1}^l\V_{{p^{\ast}}^{c_j}})$, where for each $1 \leq d_1 < d_2 <...< d_k=d$, $1 \leq c_1 < c_2 <...<c_l=d$, $\V_{p^{d_i}}$, resp. $\V_{{{p^{c_i}}^{\ast}}}$, is  $T$-invariant and is free over the algebra $\F[x]/(p(x)^{d_i})$, resp. $\F[x]/(p^{\ast}(x)^{c_i})$. We claim that $k=l$ and for each $i$, $c_i=d_i$. Using arguments similar as above, it is easy to see that $B$ is non-degenerate on $\V_{p^d} +  \V_{{p^{\ast}}^d}$. Now by induction the claim follows. Since the minimal polynomial of 
$T|_{\V_{p^{d_i}}+\V_{{p^{\ast}}^{d_i}}}$ is $p(x)^{d_i}{p^{\ast}}^{d_i}$. Hence for each elementary divisor $p(x)^{d_i}$ there is a dual elementary divisor $p^{\ast}(x)^{d_i}$ with the same multiplicity.

 Suppose $(x -1)^k$, resp. $(x+1)^k$, is an elementary divisor of $T$. Let $\W$ be a $T$-indecomposable subspace of $\V$ such that $T|_{\W}$ has characteristic polynomial  $(x-1)^k$, resp. $(x+1)^k$.  Suppose $B$ is symmetric (resp. skew-symmetric). Then it follows from \lemref{ind} that $B|_{\W}$ is either zero or non-degenerate. From \thmref{uniod} we have that $B|_{\W}$ is nondegenerate if and only if $k$ is odd (resp. even). If $B|_{\W}=0$, then $k$ must be even (resp. odd). Thus if $k$ is even (resp. odd), then by the non-degeneracy of $B$ it follows that $\W$ is a summand of a standard symmetric (resp. skew-symmetric) space, and hence the multiplicity of $(x-1)^k$ must be even (resp. odd). 

\emph{Conversely}, suppose $(i)$ and $(ii)$ of either of the theorems hold. For an elementary divisor $g(x)$, let $\V_g$ denote the $T$-indecomposable subspace  isomorphic to $\F[x]/(g(x))$.
From the theory of elementary divisors it follows that $\V$ has a decomposition 
\begin{equation}\label{x}
\V=\Sigma_{i=1}^{m_1}\V_{f_i} + \Sigma_{j=1}^{m_2} (\V_{g_i} + \V_{g_i^{\ast}}),
\end{equation} 
where  for each $i=1,2,...,m_1$,   $f_i(x)$ is either self-dual, or one of $(x+1)^k$ and $(x-1)^k$, for each $j=1,2,...,m_2$, $g_i(x)$, $g_i^{\ast}(x)$ are dual to each other and $g_i(x)\neq g_i^{\ast}(x)$.  To prove the theorem it is sufficient to induce a $T$-invariant quadratic (resp. skew-symmetric) form on each of the summands. 

Suppose $\W$ is an indecomposable summand in the above decomposition and $p(x)^k$ be the corresponding elementary divisor.

(a) Suppose $p(x)^k$ is self-dual. It follows from \propref{p-inv1} that there exists a $T$-invariant non-degenerate symmetric, as well as skew-symmetric bilinear form on $\W$. 

  Suppose $p(x)^k$ is not self-dual. Then there is a dual elementary divisor $p^{\ast}(x)^k$. As explained in \secref{standard}, there exists a $T$-invariant standard form on $\W_p + \W_{p^{\ast}}$, where $\W_p=\hbox{ker }p(T)^k$, $\W_{p^{\ast}}=\hbox{ker }p^{\ast}(T)^k$. 

(c) Suppose, $p(x)^k=(x-1)^k$. Suppose $k$ is odd (resp. even). Then the respective symmetric and skew-symmetric bilinear form is obtained from \lemref{cu} and \lemref{ind}.

Let $k=2m$, resp. $2m+1$, and the multiplicity of $(x-1)^{2m}$ is an even number. Then the number of indecomposable summands, each isomorphic to $\F[x]/(x-1)^k$, in the decomposition \eqnref{x} is even. We can pair those summands, taking two at a time, and induce a $T$-invariant standard form on each pair. 

(d) Suppose $p(x)^k=(x+1)^k$. Let $T_w$ denote the restriction of $T$ to $\W$. Then the minimal polynomial of $T_w$ is $(x+1)^k$. Thus the minimal polynomial of $-T_w$ is $(x-1)^k$. Further $T_w$ preserves a symmetric (resp. skew-symmetric) form $B$ if and only if $-T_w$ also preserves $B$. Thus this case reduces to the case (c) above, and there exists a $T$-invariant non-degenerate symmetric (resp. skew-symmetric) bilinear form on $\W$. 

This completes the proofs of \thmref{ansq1} and \thmref{ansq1'}.

\subsection{Proof of \corref{cor1}}
The proof follows from the fact that $\chi_T(x)$ is self-dual, hence $\V$ has the decomposition \eqnref{x}. Now on each of the $T$-invariant component one can induce a $T$-invariant non-degenerate symmetric or skew-symmetric form according to \thmref{ansq1} and \thmref{ansq1'}. This completes the proof.

\section{Proofs of \thmref{level} and \thmref{level2}}
Suppose $T: \V \to \V$ is unipotent. Let $\W$ be a $T$-indecomposable subspace of $\V$ of maximal dimension,  i.e. $\W$ is isomorphic to $\F[x]/(x-1)^k$.  If $B|_\W=0$, then $k \leq l$. If $k>l$, then $B|_\W$ is non-degenerate and hence $k$ must be odd  by \lemref{cu}.

Suppose $k>l$. Let $k=2m+1 \geq 3$.  Let $\W_1=\hbox{ker}(T|_{\W}-I)$. Since $\W$ is $T$-indecomposable, we must have $B|_{\W_1}=0$. Hence $r=\hbox{ dim }\W_1$ can be at most $l$. Now observe that the non-degeneracy of $\W$ implies that $k$ is at least $2r+1$. By the indecomposibility of $\W$, $k=2r+1$.  Assume $\W_1$ has the maximal dimension. Hence if the dimension of $\V$ is $2l$, then $r=l-1$ and $k \leq 2l-1$.  Suppose $\dim \V \geq 2l+1$. Then $r=l$ and $\dim \W =2l+1$. 

This completes the proof of \thmref{level}. 

The proof of \thmref{level2} is similar.

\section{The infinitesimal version: Proofs of \thmref{iinv} and \thmref{iinv1}}
We note the infinitesimal versions of the key lemmas which are crucial to the proof.  First we note the infinitesimal versions of \lemref{cu} and \lemref{ind}.

\begin{lemma}\label{icu}
Let $\V$ be a vector space over  a field $\F$ of large characteristic. Let $\V$ be equipped with a symmetric, resp. skew-symmetric, bilinear  form $B$.  Suppose $T: \V \to \V$ be a nilpotent map which keeps $B$ infinitesimally invariant.  Let $\W$ be a $T$-indecomposable subspace of $\V$.  

(i) Then  $B|_{\W}$ is either non-degenerate or $B|_{\W}=0$. 

(ii) If  $B|_{\W}$ is non-degenerate, then  the dimension of $\W$ must be odd, resp. even. 
\end{lemma}
\begin{lemma}\label{iuind}
Let $\V$ be a vector space of dimension $n\geq 2$ over a field $\F$ of large characteristic. Let $T: \V \to \V$ be a nilpotent linear map. Let $\V$ be $T$-indecomposable. Then there exists an infinitesimally  $T$-invariant non-degenerate symmetric, resp. skew-symmetric bilinear form on $\V$ if and only if $n$ is odd, resp. even. 
\end{lemma}

Next we have the following infinitesimal version of \propref{qsf}.
\begin{proposition}\label{qsfi}
Let $T$ be a linear map. If $T$ has no eigenvalue $0$, then there exists an infinitesimally $T$-invariant non-degenerate symmetric bilinear form on $\V$ if and only if there exists an infinitesimally $T$-invariant non-degenerate skew-symmetric bilinear form on $\V$. 
\end{proposition}
\begin{proof}
Suppose $B$ is an infinitesimally $T$-invariant non-degenerate symmetric bilinear form:

 for $x \ y$ in $\V$,
\begin{equation}\label{i2}B(Tx,y)+B(x, Ty)=0.\end{equation}
For $x,\ y$ in $\V$, define $B^T(x, y)=B(Tx, y)$. Then $B^T$ is a non-degenerate skew-symmetric bilinear form and it is also infinitesimally $T$-invariant:
\begin{eqnarray*}
B^T(Tx, y) +B^T(x, Ty)& =& B(T^2x, y)+B(Tx, Ty)\\
&=& -B(Tx, Ty)+B(Tx, Ty),  \ \hbox{ by  \eqnref{i2}}, \ B(T(Tx), y)=-B(Tx, Ty)\\
&=&0.
\end{eqnarray*}

The converse part follows by reversing the steps. 
\end{proof}

\begin{lemma}\label{p-inv-iv}
Let $\V$ be a vector space over a field $\F$ of large characteristic. 
Let \hbox{$S: \V \to \V$} be a linear map with characteristic polynomial $\chi_S(x)=p(x)^d$, where $p(x)$ is irreducible over $\F$ and is even.  Let $\V$ be $S$-indecomposable. Then $\dim \V$ is even, and there exists a  infinitesimaly $S$-invariant non-degenerate symmetric, as well as skew-symmetric bilinear  form on $\V$.
\end{lemma} 
\begin{proof}
Without loss of generality assume that $\V=\F[x]/((p(x)^d)$.  
Let $y=x^2$. Replacing $y$ by $x^2$ in $p(x)$, we see $p(\sqrt y)$ is a polynomial in indeterminate $y$ over $\F$, and we denote it by $q(y)$. 
Let $\E=\F[x]/(p(x))$, and $\E_1=\F[y]/(q(y))$. Clearly $\E_1$ may be taken as a subfield of $\E$, and as a field extension $\E$ has degree $2$ over $\E_1$. Since $p'(x)\neq 0$, $\E_1/\F$, and $\E/\E_1$ are separable extensions. Note that $[\E:\E_1]=2$, hence $N_{\E/\E_1}$ defines an $\E_1$-valued non-degenerate symmetric bilinear form $\n_{\E/\E_1}$ on $\E$:

for $\alpha, \ \beta \in \E_1$,
$$\n_{\E/\E_1}(\alpha,\beta)=\frac{1}{2}[N_{\E/\E_1}(\alpha+\beta)-N_{\E/\E_1}(\alpha)-N_{\E/\E_1}(\beta)].$$
It is easy to check that it is infinitesimally $S$-invariant. Now, define $B:\E \times \E \to \F$ as follows.
$$\hbox{For }\alpha, \beta \in \E, \;B(\alpha, \beta)=Tr_{\E_1/\F}(\n_{\E/\E_1}(\alpha,1)\n_{\E/\E_1}(\beta,1)).$$
Then $B$ is a non-degenerate, symmetric and infinitesimally $S$-invariant. Correspondingly, by \propref{qsfi}, there exists an infinitesimally $S$-invariant non-degenerate skew-symmetric bilinear form. 
This proves the theorem for $d=1$. 

Suppose $d$ is at least $2$. It follows from Kulkarni \cite{kulkarni} that there exists a basis such that $S=S_s+S_n$, where $S_s$ is semisimple, $S_n$ is nilpotent, $S_sS_n=S_nS_s$ and the decomposition is unique. Now the proof is similar to that of \lemref{p-inv}. 
\end{proof}
The rest of the proofs of the infinitesimal versions are similar to those of \thmref{ansq1} and \thmref{ansq1'}. We omit the details.

\end{document}